\documentclass{owrart}

\newtheorem{theorem}[equation]{Theorem}
\newtheorem{example}[equation]{Example}
\newtheorem{definition}[equation]{Definition}

\newcommand{\Tan}{\operatorname{Tan}}

\newcommand{\R}{\mathbb R}
\newcommand{\G}{{\mathbb G}}
\newcommand{\N}{\mathbb N}
\newcommand{\g}{\mathfrak{g}}
\newcommand{\Length}{\operatorname{Length}}
\newcommand{\norm}[1]{\left\Vert#1\right\Vert}

\begin{document}


\begin{talk}{Enrico Le Donne}
{Geodesic metric spaces with unique blow-up almost everywhere: properties and examples}
{Le Donne, Enrico}

\noindent

In this report we deal with metric spaces that at almost every point admit a tangent metric space. 
 These spaces are in some sense generalizations of Riemannian manifolds. We will see that, at least at the level of the tangents, there is some resemblance of a differentiable structure and of (sub)Riemannian geometry.
 I will present some results and give examples.
 
 Let $X=(X,d_X)$ and $Y=(Y,d_Y)$ be metric spaces. Fix $x_0\in X$ and $y_0\in Y$.
 If there exists $\lambda_j\to\infty$ such that,   in the Gromov-Hausdorff convergence,
 $$(X,\lambda_j d_X, x_0)\to (Y, d_Y, y_0), \qquad {\rm as \;} j\to\infty,$$
then $(Y,y_0)$ is called {\em a tangent} (or a {\em weak tangent}, or a {\em blow-up}) of $X$ at $x_0$.

Some remarks are due. Fixed $x_0\in X$, there might be more than one tangent. Moreover, in general there might not exist any tangent.
However, if the distance is doubling, then, by the work of Gromov \cite{Gromov-polygrowth}, then tangents exists. Namely, for any sequence  $\lambda_j\to\infty$, there exists 
a subsequence  $\lambda_{j_k}\to\infty$ such that $(X,\lambda_{j_k} d_X, x_0)$ converges
as $k\to \infty$. A tangent is well defined up to pointed isometry. Thus we define the set of all tangents 
of $X$ at $x_0$ as
$$\Tan(X,x_0):=\{ {\rm tangents \; of \;} X {\rm\; at\; } x_0\}/\text{pointed isometric equivalence}.$$

We consider two questions: how big is $\Tan(X,x_0)$? what happens when the tangent is unique?
The rough answer that we will give are the following. Under some `standard' assumptions,  if $(Y,y_0)\in \Tan(X,x_0)$, then
$(Y,y)\in\Tan(X,x)$, for all $y\in Y$. Moreover, in the case of unique tangents, such tangents are very special, however, not much can be said about the initial space $X$.

  \subsection*{Definition and examples}
 
 Let  $(X_j, x_j), (Y, y)$  be pointed geodesic metric spaces. We write
 $(X_j, x_j)\to (Y, y)$   in the Gromov-Hausdorff convergence if, for all $R>0$, we have
 $d_{GH}(B(x_j,R),B(y,R))\to 0.$ Here $$d_{GH}(A,B):=\inf\{d^Z_H(A',B'): Z {\rm \; metric\; space},A',B'\subseteq Z, A\stackrel{isom}{=}A', B\stackrel{isom}{=}B'\},$$
 and $d^Z_H(\cdot, \cdot)$ is the Hausdorff distance in the space $Z$.

\begin{example}
When $\R^n$ is endowed with the Euclidean distance (or more generally a norm), we have
$\Tan(\R^n, p)=\{(\R^n,0)\}, \forall p\in \R^n.$
\end{example}
\begin{example}
Let $(M,d)$ be a Riemannian manifold    (or more generally a Finsler manifold), we have
$\Tan(M,d, p)=\{(\R^n,\norm{\cdot},0)\}, \forall p\in \R^n.$
\end{example}
 
\begin{definition}[Carnot group]
Let  $\g$ be a stratified Lie algebra, i.e., 
$\g= V_1\oplus \cdots\oplus V_s, $
with $[V_j, V_1] = V_{j+1}$, for $1\leq j\leq s$, where $V_{s+1}= \{0\}$. 
 Let $\G$ be the  simply-connected Lie
group whose Lie algebra is $\g$.
Fix $\norm{\cdot}$ on $V_1$.
Define,  for any $x,y\in \G$,
\begin{equation*}\label{dist_CC}
d_{CC}(x,y):=\inf\left\{\int_0^1 {\norm{\dot\gamma(t)}}dt\;|\;\gamma\in C^\infty([0,1];\G), \gamma(0)= x , \gamma(1)= y, \dot \gamma\in V_1\right\}.
\end{equation*}
The pair $(\G, d_{CC})$ is called {\em Carnot group}.
\end{definition}

In particular, any Carnot group $\G$ is  a metric space homeomorphic to the Lie group $\G$.
Moreover, by the work of Pansu and Gromov \cite{Pansu-croissance}, the Carnot groups are the blow-downs of  left-invariant  Riemannnian/Finsler distances on $\G$. Namely, if 
$\norm{\cdot}$ is a norm on Lie$(\G)$ extending the one on $V_1$ and 
$d_{\norm{\cdot}}$ is the corresponding Finsler distance, 
$$(\G, \lambda d_{\norm{\cdot}},1)\stackrel{\lambda\to0}{\longrightarrow}(\G, d_{CC},1).$$

\begin{example} If $(\G, d_{CC})$ is a  Carnot group,
then $\Tan(\G, d_{CC},1)=\{(\G, d_{CC},1)\}.$ Indeed, for all $\lambda>0$, there is a group homomorphism $\delta_\lambda:\G \to \G$ such that $(\delta_\lambda)_*|_{V_1}$ is the multiplication by $\lambda$. Consequently, $(\delta_\lambda)_*d_{CC}=  \lambda d_{CC}.$ QED
\end{example}

  \subsection*{Results} 
Our main theorem is the following.
 \begin{theorem}[\cite{LeDonne6}] \label{main thm}
Let $(X, d)$ be a geodesic metric space. 
Let  $\mu$ be a  doubling measure. 
Assume that, for $\mu$-almost every $x\in X$, the set $\Tan(X,x)$ contains only one element.
 Then,  for $\mu$-almost every $x\in X$, the element in $\Tan(X,x)$ is a Carnot group.
 \end{theorem}
 
 \begin{example}[SubRiemannian manifolds]
 Let $M$ be a Riemannian manifold (or more generally Finsler).
 Let $\Delta\subseteq TM$ be a smooth sub-bundle. 
 Let $\mathcal X^1(\Delta)$ be the vector fields tangent to $\Delta$.
 By induction, define $\mathcal X^{k+1}(\Delta):=  \mathcal X^k(\Delta)+[\mathcal X^1(\Delta),\mathcal X^k(\Delta)]$.
 Assume that there exists $s\in \N$ such that $\mathcal X^{s}(\Delta)=TM$ and that, for all $k$, the function 
 $p\mapsto \dim \mathcal X^{k}(\Delta)(p)$ is constant.
 Define,  for any $x,y\in M$,
\begin{equation*} 
d_{CC}(x,y):=\inf\{\Length (\gamma)\;|\;\gamma\in C^\infty([0,1];M), \gamma(0)= x , \gamma(1)= y, \dot \gamma\in \Delta\}.
\end{equation*}
 Then $(M,d_{CC})$ is called an {\em (equiregular) subFinsler manifold}. 
 In such a case, by a theorem of Mitchell, see \cite{Mitchell, Margulis-Mostow},
  $$\Tan(M,d_{CC}, p)=\{ (\G, d_{CC},1)\}, \qquad \forall p\in M,$$
 with $(\G, d_{CC})$ a Carnot group, which might depend on $p$.
 \end{example}

  Theorem \ref{main thm} is proved using  the following general property.
 \begin{theorem}[\cite{LeDonne6}]\label{tantan}
Let $(X,\mu, d)$ be a doubling-measured metric space. Then, for $\mu$-almost every $x\in X$, if $(Y,y)\in \Tan(X,x)$, then
$(Y,y')\in\Tan(X,x)$, for all $y'\in Y$.
 \end{theorem}
 If $\#\Tan(X,x_0)=1$, then $(Y, y_0)=(Y,y),$ for all $y\in Y$. In other words,  the isometry group  Isom$(Y)$ acts on $Y$ transitively. Thus we use the following.
 
 \begin{theorem}[Gleason-Montgomery-Zippin, \cite{mz}]\label{Montgomery-Zippin}
Let $Y$  be a metric space that is complete, proper, connected, and locally connected.
Assume that the isometry group Isom$(Y)$ of $Y$ acts transitively on $Y$.
Then Isom$(Y)$  is a Lie group with finitely many
connected components.
\end{theorem}
Regarding the conclusion of the proof of   Theorem \ref{main thm},
 since moreover $Y$  is geodesic, being $X$ so, then $Y$ is  a subFinsler manifold,  by \cite{b}.
From Mitchell's Theorem  and the fact that $\{Y\}=\Tan(Y,y)$, $Y$ is a Carnot group. QED

\subsubsection*{Comments and more examples} There are other settings in which the tangents are (almost everywhere) unique.
The snow flake metrics $(\R, \norm{\cdot}^\alpha)$ with $\alpha\in (0,1)$ are such examples.
Some examples on which the tangents are Euclidean spaces are the
Reifenberg vanishing flat metric spaces, which have been considered in 
\cite{Cheeger-Colding,David-Toro}. Alexandrov spaces have Euclidean tangents almost everywhere, \cite{Burago-Gromov-Perelman}.

However, even in the subRiemannian setting, the tangents are not local model for the space. Indeed, there are subRiemannian manifolds  with a different tangent at each point,
\cite{Varchenko}.
In fact, 
there exists a nilpotent Lie group equipped with left invariant sub-Riemannian metric that is not locally biLipschitz equivalent to its tangent, see \cite{LOW}.
Such last fact can be seen as the local counterpart of 
a result by Shalom, which states that  
there exist  two  finitely generated nilpotent   groups $\Gamma$ and $\Lambda$ that have the same blow-down space, but they are not   quasi-isometric equivalent, see  \cite{Shalom}.

Another pathological example from \cite{Hanson-Heinonen} is the following.
For any $n>1$, there exists a geodesic space $X$ supporting a doubling measure $\mu$
such that at $\mu$-almost all point of $X$ the tangent is $\R^n$, but $X$
has no manifold points.

\end{talk}

\end{document}